\documentclass{amsart}
\usepackage{epsfig}
\usepackage{amsxtra,amssymb,amsthm,amsmath,latexsym}

\newtheorem{theorem}{Proposition}[section]
\newtheorem{lemma}[theorem]{Lemma}

\theoremstyle{definition}

\newtheorem{fact}[theorem]{Fact}

\theoremstyle{remark}
\newtheorem{remark}[theorem]{Remark}

\numberwithin{equation}{section}

\def\qed{\hfill\ \vbox{\hrule\hbox{\vrule\kern4pt\vbox{\kern4pt{}
\kern4pt}\kern4pt\vrule}\hrule}}
\def\cirN{\buildrel\circ\over N}
\def\disju{\ {\perp\!\perp}\ }

\begin{document}

\title{Topologically knoted Lagrangians
in simply connected four manifolds
}

\author{Dave Auckly}

\address{Department of Mathematics,\\ 
Kansas State University,\\
Manhattan, KS 66506-2602, USA}

\thanks{This work was partially supported by National Science Foundation 
grant DMS-0204651.}

\date{June ??, 2003 and, in revised form, ??.}

\begin{abstract}

Vidussi was the first to construct  knotted Lagrangian tori
in simply connected four dimensional manifolds. Fintushel and Stern introduced a second way to detect such knotting.
This note demonstrates that similar examples may be  distinguished by the fundamental group of the exterior. 
\end{abstract}

\maketitle

\section*{Introduction}
The study of the topology of Lagrangian submanifolds took off starting in the mid 1980's. This early work is surveyed in an influential paper of Eliashberg and Polterovich \cite{EP}. The existence of smoothly knotted Lagrangian tori in 
some simply connected four dimensional symplectic manifolds was recently
demonstrated by S. Vidussi, \cite{V}. The paper of Vidussi contains more history
and background of the problem. Shortly thereafter, R. Fintushel and R. Stern
introduced a second invariant that could distinguish smoothly knotted Lagrangian
tori and applied their invariant to variations of Vidussi's construction \cite{FS}.
In this paper we note that the fundamental group of the complement of a Lagrangian
torus may also be used to show that it is knotted. 

I would like to thank  Stefano Vidussi for asking the question that led to this paper,
and for several discussions about this material. I would also like to thank Ron Fintushel for his comments on a preliminary version of the paper.

\section*{Construction of knotted Lagrangian tori and the invariants of Vidussi, and Fintushel and Stern}

The essence of the construction of knotted Lagrangians that was devloped by Vidussi
in \cite{V} may be summarized as follows. Let $M$ be a $3$-manifold that fibers
over $S^1$. The $4$-manifold $S^1\times M$ inherits a symplectic form expressible as $\omega= \pi_1^*d\theta\wedge p^*d\theta+\eta$, where $d\theta$ is the standard $1$-form
on $S^1$, $\pi_1$ is the natural projection $S^1\times M\to S^1$, $p$ is the composition of the projection onto the second factor followed by the projection of $M$ onto $S^1$, and  $\eta$ is an area form on the fibers of $M$.  
In this manifold, any torus of the form, $S^1\times \delta$ is symplectic i.e. $\omega$
restricts to a symplectic form whenever $\delta$ is an embedded circle transverse to all of the fibers of $M$. Furthermore, any torus of the form, $S^1\times \gamma$ is Lagrangian i.e. $\omega$
restricts to zero whenever $\gamma$ is an embedded circle in one of the fibers of $M$. More complicated examples may be constructed by taking symplectic fiber sums with other symplectic manifolds. 

In \cite{V} Stefano Vidussi  used the Seiberg-witten invariants of the symplectic sum of an $E(1)$ with those examples
summed with non-rational elliptic surfaces
along  the Lagrangian torus to prove that the various
Lagrangian tori were not smoothly isotopic. Ron Fintushel and Ron Stern later used
relative Seiberg Witten invariants of the complement of the Lagrangian tori to prove that they are not smoothly isotopic \cite{FS}. Fintushel and Stern described a particularly
simple special case of their invariant as a Lagrangian framing defect. 

\section*{Fundamental groups of Lagrangian knot complements}

The fundamental group of the complement of the torus is a natural invariant that
could in principal detect (even topologically) non-isotopic  tori. For many examples
the fundamental group of the complement will not detect anything. However there 
are homotopic Lagrangian tori having complements with non isomorphic fundamental groups. These examples may be chosen to be null-homologous with zero Lagrangian
framing defect.

The first specific example considered in \cite{V} took the following form. Let $T(p,q)$ represent the $(p,q)$ torus knot, and $M_{\tau}$ be the result of $0$-surgery on $\tau=T(2,-3)$.
It is well known that $M_{\tau}$ fibers over $S^1$ \cite{R}. Let $E(1)$ be the 
rational elliptic surface \cite{GS}, and set $E(1)_\tau=E(1)\#_{F=S^1\times \mu}S^1\times M_{\tau}$, where $\mu$ is a meridian to $\tau$. Figure 1 below shows, among other things, a knot $\gamma_p$ on a fiber in the exterior of the knot $\tau$. For the first example, we will ensure that $\gamma_p$ does not link $\mu$. By inspection, $\gamma_p = T(p,p+1)$ as a knot in $S^3$. As we will see the  fundamental group of the complement of $S^1\times \gamma_p$ in $E(1)_\tau$ is isomorphic to ${\mathbb Z}$.
However  doubling this construction in an appropriate sense will result in a fundamental group that surjects onto the fundamental group of the complement
of the $(p,p+1)$ torus knot. 
\bigskip

\hskip 25bp \epsfig{file=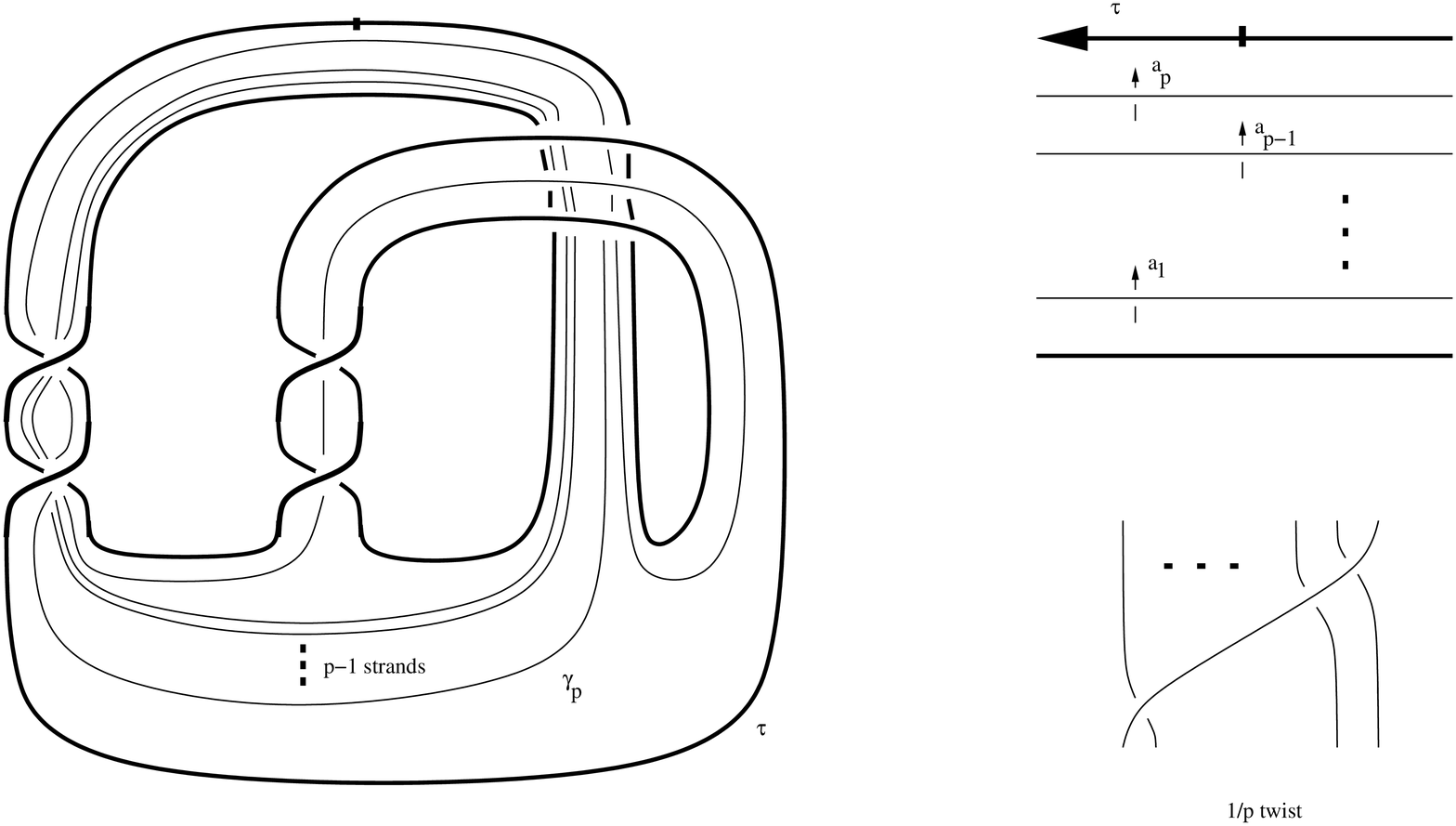 , width=4truein}

\smallskip
\centerline{Figure 1: The knot $\gamma_p$}
\smallskip

Our doubled construction begins with the knot, $\tau\#-\overline\tau$. Here the overline refers to taking the reflection of the knot, and the minus sign refers to
switching the orientation of the knot. Our symplectic manifold is defined by,
$E(1)_{\tau\#-\overline\tau}=E(1)\#_{F=S^1\times \mu}S^1\times M_{\tau\#-\overline\tau}$, where $\mu$ is a meridian to $\tau\#-\overline\tau$. We will
consider  the Lagrangian tori, $S^1\times (\gamma_p\#-\overline\gamma_p)$
in $E(1)_{\tau\#-\overline\tau}$.

\bigskip

\hskip 25bp \epsfig{file=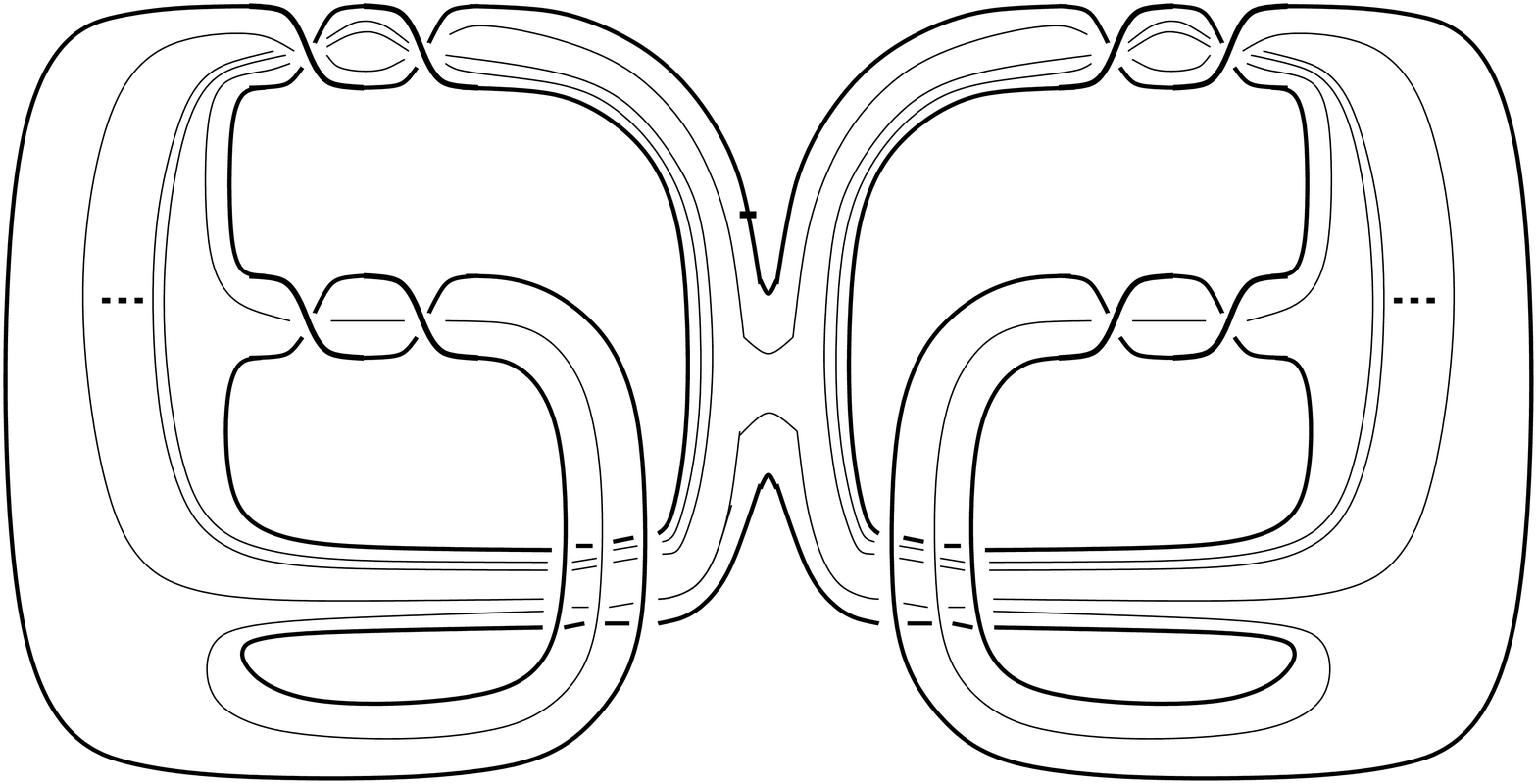 , width=4truein}

\medskip
\centerline{Figure 2: The double, $\tau\#-\overline\tau\disju\gamma\#-\overline\gamma$}
\smallskip
\begin{theorem}\label{2gp}
The fundamental groups of the Lagrangian knot exteriors $S^1\times (\gamma_p\#-\overline\gamma_p)$ are given by:
$$
\begin{aligned}
\pi_1(E(1)_{\tau\#-\overline\tau}&-{\buildrel\circ\over N}(S^1\times   (\gamma_p\#-\overline\gamma_p)))\\
&=\langle u, v, x, y | u^pv^{p+1}=1, x^py^{p+1}=1, uv=xy, vu=yx \rangle.
\end{aligned}
$$ 
\end{theorem}

\begin{proof} We will use the following alternate description of $E(1)_{\tau\#-\overline\tau}$. Let $X_{\tau\#-\overline\tau}$ denote the exterior of $\tau\#-\overline\tau$. We then have $E(1)_{\tau\#-\overline\tau}=S^1\times X_{\tau\#-\overline\tau}\cup_{T^3}(E(1)-\cirN(T^2))$. Here $T^2$ is the fiber in $E(1)$. This allows
us to apply the Seifert-Van Kampen theorem to compute $\pi_1(E(1)_{\tau\#-\overline\tau}-\cirN(S^1\times\gamma_p\#-\overline\gamma_p))$ as follows.
$$
\begin{aligned}
\pi_1 & (E(1)_{\tau\#-\overline\tau}-\cirN(S^1\times\gamma_p\#-\overline\gamma_p)) \\
& =\pi_1(S^1\times (S^3-\cirN(\tau \disju\gamma_p\#-\overline\gamma_p))\cup_{T^3}(E(1)-\cirN(T^2))) \\
&= \pi_1(S^1\times (S^3-\cirN(\tau\#-\overline\tau \disju\gamma_p\#-\overline\gamma_p)))/\langle\pi_1{T^3}\rangle \\
&=\pi_1(S^1\times ((S^3-\cirN(\tau\#-\overline\tau \disju\gamma_p\#-\overline\gamma_p))\cup_{T^2}S^1\times D^2))
\cup_{\buildrel{S^1\times pt}\over{ pt\times\tau\#-\overline\tau }}2D^2) \\
&=\pi_1((S^3-\cirN(\gamma_p\#-\overline\gamma_p))\cup_{\tau\#-\overline\tau}D^2) = \pi_1(S^3-\cirN(\gamma_p\#-\overline\gamma_p))/\langle\tau\#-\overline\tau\rangle. 
\end{aligned}
$$ 
The third line of the above computation follows because $\pi_1(E(1)-\cirN(T^2))=1$. The fourth line follows because glueing anything that kills the fundamental group of the boundary will have the same effect. The $(S^1\times D^2)\cup D^2 \cup D^2$ that we glue in replaces $\cirN(\tau\#-\overline\tau)$, kills the generator
of the $S^1$ factor, and adds one more relation. To complete the proof we need to
compute $\tau\#-\overline\tau$  in $\pi_1(S^3-\cirN(\gamma_p\#-\overline\gamma_p))$.

We begin by computing the element of $\pi_1(S^3-\cirN(T(p,p+1)))$ represented by $\tau$. We will use the 
Wirtinger generators depicted in figure 1.
Recall that a Wirtinger generator is a loop that starts at the base point, goes over 
the knot, under in the direction of the arrow and then back over the knot to the base point, \cite{R}. Using the Wirtinger relation we see that the generators on the bottom
of a $1/p$ twist will be $a_1a_1a_1^{-1}$, $a_1a_pa_1^{-1}$, \dots $a_1a_2a_1^{-1}$, when the generators at the top of the twist are $a_p$, \dots $a_1$. Since $T(p,p+1)$ is just a stack of $p+1$ $1/p$-twists, this allows us to obtain the presentation, $\langle z, a_1,\dots a_p | z=a_1\dots a_pa_1, za_1z^{-1}=a_p, za_2z^{-1}=a_1, \dots za_pz^{-1}=a_{p-1}\rangle$, and label every strand in the diagram for
$T(p,p+1)$. The element of $\pi_1$ represented by $\tau$ may now be read off of the 
diagram. We obtain, $\tau= (y^{-1})(ya_1y^{-1})y(ya_1y^{-1})(ya_1^{-1}y^{-1})(ya_1^{-1}y^{-1})=[a_1,y]$ where $y=a_1\dots a_p$. In computing this expression, we reversed the obvious Reidemeister I move on the right side of
figure 1, and grouped terms after each group of undercrossings. Since $\gamma_p$
lies on a Seifert surface for $\tau$, we know that the linking number of $\gamma_p$
and $\tau$ is zero, so $\tau$ must be (as we found it) in the commutator subgroup of the knot group of
$\gamma_p$. In this case $\tau$ normally generates the commutator subgroup, so
the knot group of $\gamma_p$ mod $\tau$ is isomorphic to the infinite cyclic group.
The standard generators of the knot group of $T(p.p+1)$
are $x=a_1^{-1}y^{-1}$ and $y$.

We will apply the Seifert VanKampen theorem to compute the knot group of
$\gamma_p\#-\overline\gamma_p$, and the element $\tau\#-\overline\tau$
in that group. The twice punctured plane of symmetry in the complement of
$\gamma_p\#-\overline\gamma_p$ in figure 2 becomes a twice punctured $S^2$
when we add the point at infinity. Two copies of the complement of $\gamma_p$
are glued together along this punctured sphere. The resulting presentation of the fundamental group is
$$
\begin{aligned}
\langle w, z, & a_k, b_k, k=1,\dots p | \\
 & z=a_1\dots a_pa_1, w=b_1\dotsb_pb_1, za_{k+1}z^{-1}=a_k, wb_{k+1}w^{-1}=b_k, k=1,\dots p\rangle,
\end{aligned}
$$
and $\tau_p\#-\overline\tau_p$ represents $[a_1,y][b_1,v]^{-1}$ where $v=b_1\dotsb_p$. Setting $x=z^{-1}$,  $u=w^{-1}$, and rewriting the 
presentation of the knot group of $\gamma_p\#-\overline\gamma_p$
in terms of $u$, $v$, $x$, and  $y$ results in the stated presentation
of the complement of $S^1\times(\gamma_p\#-\overline\gamma_p)$.
\end{proof}

\begin{remark}
The computation at the beginning of the previous lemma shows that
the fundamental group of the complement of a Lagrangian torus of the form
$S^1\times\gamma$ in $E(1)_\kappa$ always takes the form $\pi_1(X_\gamma)/\langle\kappa\rangle$ where $\kappa\in [\pi_1(X_\gamma),\pi_1(X_\gamma)]$ is the class of $\kappa$ and $\gamma$
is a knot in a fiber of the fibered knot, $\kappa$.
\end{remark}

\begin{remark}
The fundamental group of the complement of $S^1\times\gamma_p$ in $E(1)_\tau$
is always ${\mathbb Z}$.  However the Fintushel-Stern Lagrangian framing defect is given by $\lambda(\gamma_p)=p+1$, so the tori $S^1\times\gamma_p$ represent distinct isotopy classes. 
\end{remark}

\begin{remark}For  any knot of the form, $\gamma\#-\overline\gamma$ in a fiber of a knot of the form $\kappa\#-\overline\kappa$, the fundamental group of the complement
of $S^1\times(\gamma\#-\overline\gamma)$ in $E(1)_{\kappa\#-\overline\kappa}$ will surject onto the knot group of $\gamma$. Thus such Lagrangian tori are not isotopic to the standard
nullhomologous Lagrangian tori. By symmetry
the Lagrangian framing defect of Fintushel and Stern will vanish on the tori
$S^1\times(\gamma\#-\overline\gamma)$. Vidussi's technique of summing with an
$E(1)$ is an appropriate way or Fintushel and Stern's main invariant $I(X,T)$ may still in principal detect this.
\end{remark}

\begin{remark}The next proposition will establish that all of the groups in the previous
proposition are distinct. Thus the tori $S^1\times(\gamma_p\#-\overline\gamma_p)$
form an infinite family of nonisotopic Lagrangian
tori that are all homotopic.
\end{remark} 

\begin{lemma}\label{dif}
The following groups are all distinct.
$$
\Gamma_p=\langle u, v, x, y | u^pv^{p+1}=1, x^py^{p+1}=1, uv=xy, vu=yx \rangle
$$
\end{lemma}

\begin{proof} We will use the order ideal of the Alexander module to distinguish these
groups \cite{R}. Introducing new variables, $t=xy$, $a=t^pv$, and $b=t^py$ will
allow the presentation to be written in the form:
$$
\begin{aligned}
\langle t, a, b | & \prod_{k=0}^{p-1} t^{k(p+1)+1}a^{-1} t^{-(k(p+1)+1)}
\prod_{k=0}^p t^{p^2-kp}a t^{-(k(p+1)+1)}, \\
& \prod_{k=0}^{p-1} t^{k(p+1)+1}b^{-1} t^{-(k(p+1)+1)} 
 \prod_{k=0}^p t^{p^2-kp}b t^{-(k(p+1)+1)}, ata^{-1}t^{-1}tbt^{-1}b^{-1} \rangle.
\end{aligned}
$$
From this presentation it is clear that the abelianization of $\Gamma_p$ is infinite
cyclic. Furthermore from the same presentation one can see that the commutator
subgroup, ${\mathcal D}^1\Gamma_p = [\Gamma_p,\Gamma_p]$ is generated by the elements, $t^kat^{-k}$,
and $t^kb t^{-k}$ for all $k$. The abelianization of the commutator subgroup,
$({\mathcal D}^1\Gamma_p)^{\hbox{ab}}$ inherits the structure of a ${\mathbb Z}[t,t^{-1}]$-module
with $t$ acting by conjugation. As a ${\mathbb Z}[t,t^{-1}]$-module it has a presentation,
$$
({\mathcal D}^1\Gamma_p)^{\hbox{ab}}=\langle a, b, | {\mathfrak p}(t)a=0, {\mathfrak p}(t)b=0, (1-t)a-(1-t)b=0 \rangle,
$$
here 
$$
{\mathfrak p}(t)=\sum_{k=0}^p t^{p^2-kp}-\sum_{k=0}^{p-1} t^{k(p+1)+1}
=\frac{(t^{p(p+1)}-1)(t-1)}{(t^{p+1}-1)(t^p-1)}.
$$
We are making the obvious abuses of notation with $t$, $a$, and $b$. This
module is called the Alexander module. The order ideal of this module
is the ideal of ${\mathbb Z}[t,t^{-1}]$ generated by ${\mathfrak p}(t)^2$ and $(t-1){\mathfrak p}(t)$.
The number $e^{\frac{2\pi i}{p(p+1)}}$ is a zero of every polynomial in
the order ideal, and the numbers $e^{\frac{2\pi i}{k(k+1)}}$ are not zeros
of ${\mathfrak p}(t)^2$ for any $k>p$. This completes the proof.
\end{proof}

\bibliographystyle{amsplain}

\end{document}